\documentclass[12pt]{amsart}
\usepackage[cp1251]{inputenc}
\usepackage[russian]{babel}
\usepackage{amssymb}

\theoremstyle{plain}

\newcommand{\intl}{\mathop{\int}\limits}

\newcommand{\suml}{\mathop{\sum}\limits}

\renewcommand{\le}{\leqslant}
\renewcommand{\ge}{\geqslant}

\def\al{\alpha}
\def\la{\lambda}

\def\eps{\varepsilon}

\def\w{{\mathrm w}}
\def\v{{\mathrm v}}

\newlength{\lenun}
\newlength{\lendu}
\makeatletter \@addtoreset{equation}{section} \makeatother

\begin{document}

\title[]{Свойства отображения, связанного с
восстановлением \\
оператора Штурма-Лиувилля по спектральной функции.\\
Равномерная устойчивость в шкале соболевских\\
пространств.}
\thanks{Работа поддержана Российским фондом фундаментальных
исследований, грант № 10-01-00423 }

\author[А.М.Савчук, A. A. Шкаликов] {А.М.Савчук, А.А.Шкаликов}


\maketitle

{\bf  Abstract}. {\it Denote by $L_D$ the Sturm-Liouville operator $Ly=-y'' +q(x)y$ on the finite interval $[0,\pi]$
with Dirichlet boundary conditions $y(0)=y(\pi)=0$. Let $\{\lambda_k\}_1^\infty$ and $\{\alpha_k\}_1^\infty$ be the
sequences of the eigenvalues and norming constants of this operator. For all $\theta \geqslant 0$ we study the map $F:
W_2^{\theta} \to l_D^\theta$ defined by $F(\sigma) =\{s_k\}_1^\infty$. Here $\sigma= \int q $ is the primitive of $q$,\
$\bold s = \{s_k\}_1^\infty$ be regularized spectral data defined by $s_{2k} =\sqrt{\lambda_k}-k,\
s_{2k-1}=\alpha_k-\pi/2$ and $l_D^\theta$ are special Hilbert spaces which are constructed in the paper as finite
dimensional extensions \ of \  the usual weighted $l_2$ spaces. We give a complete characterization of the image of
\ this nonlinear operator, show that it is locally
invertible analytic map, find explicit form of its
Frechet derivative. The main result of the paper are
the uniform estimates of the form
$\|\sigma-\sigma_1\|_\theta
\asymp \|\bold s -\bold s_1\|_\theta$, provided that
the spectral data $\bold s$ and $\bold s_1$ run
through special convex sets in the spaces
$l_D^\theta$.}

\bigskip
Настоящая работа связана с изучением обратной задачи
для оператора Штурма--Лиувилля с краевыми условиями
Дирихле:
\begin{equation}\label{1}
Ly=-y''+q(x)y,\qquad x\in[0,\pi],\qquad
y(0)=y(\pi)=0.
\end{equation}
Решение задачи о восстановлении потенциала $q$  по спектральной функции  этого оператора (далее обозначаем его через
$L_D$)  было дано в классической работе Гельфанда и Левитана
\cite{GL}, которой предшествовали классические работы
Борга \cite{Bo}  и Марченко \cite{Ma1}. Впоследствии
появились сотни работ на эту тему, обзор работ и
ссылки можно найти в книге Фрайлинга и Юрко
\cite{FrYu}, а с учетом работ последнего десятилетия, ---
в недавней работе авторов \cite{SS5}.

В работе авторов \cite{SS1} было предложено определение операторов Штурма-Лиувилля для потенциалов--распределений $q$
из пространства Соболева $W^{-1}_2[0,\pi]$  и было предпринято изучение спектральных характеристик таких операторов.
Решение рассматриваемой обратной задачи для сингулярного случая было дано Гринивым и Микитюком \cite{HM2}.{\it Цель
этой статьи
--- провести анализ прямой и обратной задач для оператора $L_D$  с потенциалами из всей шкалы пространств Соболева
$W^{\alpha}_2[0,\pi], \ \,
\alpha \geqslant -1$ и доказать равномерные двусторонние оценки
для  разности потенциалов в норме пространства $W^{\alpha}_2$ (при $\al>-1$) через норму спектральных данных в
конструируемых нами пространствах $l_D^{\,\alpha +1}$, куда помещаются регуляризованные спектральные данные.} Для этой
цели мы изучаем свойства нелинейного отображения $F:\, W^{\alpha +1}_2 \to l_D^{\,\alpha +1}$, которое ставит в
соответствие первообразной потенциала $\sigma =\int q$ регуляризованные спектральные данные оператора $L_D$. В нашей
предыдущей работе \cite{SS4} мы уже начали изучение этого отображения, показав, что при $\alpha
> -1$ оно является слабо нелинейным (т.е. компактным
возмущением линейного отображения). Здесь мы получим существенно более полную информацию об этом отображении. В
частности, мы даем точное описание (характеризацию) спектральных данных, когда потенциал $q$ пробегает $
W^{\alpha}_{2,\mathbb R}$ (так мы обозначаем  все вещественные функции из пространства $ W^{\alpha}_2$),  приводим
явные формулы для восстановления первообразной потенциала (леммы 2.1 и 2.2),  находим явный вид производных прямого и
обратного отображений. {\it Основные результаты работы сформулированы в теоремах 2.6, 2.8, 2.10 и 2.15.}

Отметим, что язык теории аналитических отображений
при исследовании обратных задач для оператора
Штурма-Лиувилля получил развитие в работах Трубовица
и его соавторов (см. подробное изложение в книге
\cite{PoTr}). Однако изучаемое нами отображение
отличается от отображений, расмотренных в
\cite{PoTr}.  \ Кроме того, изучаем мы его
одновременно во всей шкале соболевских пространств с
ипользованием теории нелинейной интерполяции. Поэтому
получаемые нами формулы другие, нежели в
\cite{PoTr} и доказательство их проводится на другом
пути.

{\it Эту статью следует рассматривать как дополнение к недавней работе авторов \cite{SS5},} в которой намеченный план
исследования был полностью реализован для отображения
$F=F_B$,  связанного с обратной задачей Борга
восстановления потенциала по двум спектрам.  Там же в
\cite{SS5} показано, что схема исследования
отображения $F_B$ полностью сохраняется при
исследовании отображения $F=F_D$, связанного с
обратной задачей восстановления потенциала по
спектральной функции. В \cite{SS5} были
сформулированы леммы и теоремы о свойствах
отображения $F_D$,  однако в виду ограничения объема
статьи доказательства основных лемм были опущены.
Поэтому {\it основной целью статьи можно считать
представление полных доказательств сформулированных в
\cite{SS5} результатов об отображении $F_D$.}
Для удобства читателя мы приводим здесь заново все
необходимые сведения об операторе $L_D$ и формулируем
ранее полученные результаты об отображении $F=F_D$.
Поэтому статью можно читать независимо от
\cite{SS5} с одной оговоркой: здесь не делается
сравнений с полученными ранее результатами и список
ссылок сокращен до минимума (к сожалению, авторам не
удалось сократить число ссылок  на свои работы).
Следует также иметь ввиду, что далее мы работаем не с
отображением $q\to
\{\text{спектральные данные}\}$, а с отображением
$\sigma =\int q \, \to\{\text{спектральные
данные}\}$. Условие $q\in  W^{\alpha}_2,\ \,
\alpha\geqslant -1$ при этом переходит в условие
$\sigma \in  W^{\,\theta}_2,\ \,
\theta\geqslant 0$. Классический случай $q\in L_2$
будет соответствовать значению $\theta= 1$. Индекс $D$ в записи отображения  $F=F_D$ будем опускать.

\section{\bf Определение пространства $l_D^{\,\theta}$ и
отображения $F$. Предварительные сведения о
 свойствах $F$.
 }

 Сначала напомним, что определение
оператора Штурма--Лиувилля  с классическим потенциалом $q\in L_1[0,\pi]$ можно расширить для потенциалов--распределений
из cоболевского пространства $W^{-1}_2[0,\pi]$. Предположим, что комплекснозначный потенциал $q$ принадлежит
соболевскому пространству $W^\alpha_2[0,\pi]$ при некотором $\alpha
\geqslant -1$. Положим $\sigma(x) = \int q(x)\, dx$,
где первообразная понимается в смысле распределений. Согласно
\cite{SS1}
 определим
оператор Дирихле равенством
\begin{equation}
\label{2}
L_Dy=Ly=-(y^{[1]})'-\sigma(x)y^{[1]}-\sigma^2(x)y,\quad
y^{[1]}(x):=y'(x)-\sigma(x)y(x),
\end{equation}
взяв в качестве области определения
$$
\mathcal D(L_D)=\{y,\, y^{[1]}\in W_1^1[0,\pi]\ \,\vert\ \,  Ly\in L_2[0,\pi],\ y(0)=y(\pi)=0\}.
$$
Для гладких функций $\sigma$ правые части в
\eqref{1} и \eqref{2} совпадают и мы получаем
классический оператор Штурма-Лиувилля с краевыми
условиями Дирихле. Отметим, что добавление к функции
$\sigma$  константы не меняет оператор $L_D$.

Обозначим через ${\rm s}(x,\lambda)$ единственное решение уравнения  $Ly-\lambda y =0$, удовлетворяющее условиям
${\rm s}(0,\lambda) = 0$ и ${\rm s}^{[1]}(0,\lambda)= \sqrt\lambda$ (известно, см. \cite{SS1}, что такое решение существует и единственно). Здесь и далее мы выбираем аргумент корня $\arg\sqrt\la\in(-\pi/2,\pi/2]$. Очевидно, нули
$\{\lambda_k\}_1^\infty$ целой функции ${\rm s}(\pi,\lambda)/\sqrt{\la}$ являются собственными значениями оператора $L_D$. В случае вещественного потенциала $q$ все собственные значения являются простыми и вещественными, в этом случае считаем их занумерованными так, чтобы последовательность $\{\lambda_k\}_1^\infty$ была строго возрастающей. Для комплексных $q$ нумерацию можно провести так, чтобы последовательность $\{|\lambda_k|\}_1^\infty$ не убывала. Введем в рассмотрение также числа
\begin{equation*}
\alpha_k =\begin{cases}\int_0^\pi {\rm s}^2(x,\lambda_k)\, dx,\quad\qquad\text{если}\ \la_k\ne0;\\
\int_0^\pi\left(\tfrac{{\rm s}(x,\lambda)}{\sqrt\la}\right)^2\, dx
\Big\vert_{\la=\la_k},\ \ \text{если}\ \la_k=0,
\end{cases}
\end{equation*}
которые в случае вещественных потенциалов называются
{\it нормировочными константами} (мы сохраним это название и для
комплексных потенциалов). Последовательности
\begin{equation*}
\{\lambda_k\}_1^\infty \cup \{\alpha_k\}_1^\infty
\end{equation*}
формируют спектральные данные оператора $L_D$. Нетрудно заметить
(см., например, \cite[гл.1]{Le}), что задание этих
данных для вещественных потенциалов эквивалентно заданию спектральной функции оператора $L_D$.

Регуляризованные спектральные данные определим
следующим образом
\begin{equation}\label{L_D}
s_{2k} = \sqrt{\lambda_k} -k,\qquad s_{2k-1} = \alpha_k -\pi/2,\qquad k=1,2, \dots.
\end{equation}

Теперь  построим пространства, которым  принадлежат
регуляризованные спектральные данные. Обозначим через
$l^{\,\theta}_2$ весовое $l_2$-пространство,
состоящее из последовательностей комплексных чисел
$\bold x=\{x_1, x_2,\dots\}$, таких, что
$$
\|\bold x\|^2_\theta : =\sum_1^\infty |x_k|^2\, k^{2\theta} <\infty.
$$
Рассмотрим специальные последовательности
\begin{gather*} \bold e_{2s-1} =
\{\, 0,\  2^{-(2s-1)},\, 0,\ 4^{-(2s-1)},\,
0,\ 6^{-(2s-1)},\ldots\}\qquad\text{и}\\ \quad
\bold e_{2s}=
\{2^{-(2s)},\, 0,\  4^{-(2s)},\, 0,\  6^{-(2s)},\ldots\},
 \qquad s=1,2,\dots .
\end{gather*}
Заметим, что последовательность $\bold e_p $ принадлежит пространству $l_2^{\,\theta}$ при $0\leqslant\theta < p-1/2$ и
$\bold e_p $ не принадлежит $l_2^{\,\theta}$ при $\theta
\geqslant p-1/2$. Для фиксированного $\theta\geqslant
0$   однозначно определяется целое число $m$,
подчиненное условию $m-1/2
\leqslant \theta < m+1/2$.  Для такого $\theta$ определим
пространство $l_D^{\,\theta}$ как конечномерное расширение
пространства $l_2^{\,\theta}$ следующим образом
$$
l_D^{\,\theta}=l_2^{\,\theta}\oplus
\text{span}\{\bold e_k\}_{k=1}^{m}.
$$
Таким образом,  $l_D^{\,\theta}$ состоит из элементов
$\bold{x}+\sum_{k=1}^{m} c_k \bold e_k$, где $\bold
x\in l_2^{\,\theta}$, а $\{c_k\}_1^{m}$ ---
произвольные комплексные числа. Скалярное
произведение элементов из $l_D^{\,\theta}$
определяется формулой
$$
(\bold{x} +\sum_{k=1}^{m} c_k \bold e_k , \ \bold{y}
+\sum_{k=1}^{m} d_k \bold e_k ) = (\bold x ,
\bold y)_\theta + \sum_{k=1}^{m}
c_k\overline{d_k},
$$
где $(\bold x , \bold y)_\theta$ --- скалярное произведение в $l_2^{\,\theta}$. Построенное пространство  свяжем с
регуляризованными спектральными данными для оператора $L_D$. Хотя это пространство определено как конечномерное
расширение весового пространства $l_2^{\,\theta}$, его элементы удобнее записывать в форме обычных последовательностей.
Например, при $3/2\leqslant
\theta <5/2$  элементы пространства $l_D^{\,\theta}$
состоят из последовательностей $\{x_k\}_1^\infty$ с
координатами
$$
x_k =y_k +\begin{cases}
\alpha_1\, (k+1)^{-2},\quad&
\text{если }k\text{ нечетно,}\\
\alpha_2 \, k^{-1},\quad& \text{если }k\text{ четно,}
\end{cases}\quad \text{где}\ \ \{y_k\}_1^\infty \in
l_2^{\,\theta}, \ \ \alpha_1, \alpha_2 \in \mathbb C.
$$
Из такого представления легко следует, что пространство $l_D^{\,\eta}$  компактно вложено в пространство
$l_D^{\,\theta}$  при $\eta > \theta$ (здесь мы принимаем во внимание компактность вложения
$l^{\,\eta}_2\hookrightarrow l_2^{\,\theta}$\ при $\eta > \theta$).

Определим  нелинейный оператор
\begin{equation}\label{defi}
F (\sigma) = \{s_k\}_1^\infty ,
\end{equation}
где числа $\{s_k\}$ определены равенствами
\eqref{L_D}.
Из результатов работ
\cite{SS2} и
\cite{HM1}
следует, что последовательности, образованные из регуляризованных спектральных данных  в правых частях равенств
\eqref{L_D} являются последовательностями из $l_2$ для любой первообразной  $\sigma =\int q(x)\, dx\in L_2(0,\pi)$.
Поэтому определенный равенством
\eqref{defi} оператор $F$ корректно определен как
оператор из $L_2$  в $l_2$. Более того,  согласно результатам \cite{SS3} и \cite{SS5}, образ сужения этого оператора на
соболевское  пространство $W^\theta_2,\ \, \theta >0,$ лежит в пространстве $l_D^{\,\theta}$. Именно для этой цели мы
проводили расширения пространств  $l_2^{\,\theta}$. Без присоединения к  $l_2^{\,\theta}$ специальных
последовательностей соответствующий результат неверен.

Важную роль в дальнейшем играют результаты,
доказанные в работе
\cite{SS4},  которые для удобства приведем
здесь в нужном нам виде.

{\bf Теорема 1.1} {\sl При любом фиксированном $\theta \geqslant 0$ нелинейный оператор $F$ корректно определен как
оператор из пространства $W^\theta_2$ в  $l_D^\theta$ и дифференцируем по Фреше в каждой точке (функции) $\sigma$ при
условии, что эта функция вещественнозначна, а числа $\la_k(\sigma)$ не равны нулю. В частности, этот оператор
дифференцируем по Фреше в точке $\sigma=0$, причем  производная по Фреше в этой точке есть линейный оператор $T$,
определяемый равенством}
\begin{equation*}
\begin{cases}
(T\,\sigma)_{2k-1}=-\int\limits_0^\pi(\pi-t)\sigma(t)\cos(2kt)\, dt ,
\quad k=1,2,\dots,\\
(T\,\sigma)_{2k}=-\tfrac1{\pi}\int\limits_0^\pi\sigma(t)\sin(2kt)\, dt,\quad k=1,2,\dots.
\end{cases}
\end{equation*}

{\bf  Доказательство}  этого утверждения получается из Предложения
1 и Теоремы 4.2 работы  \cite {SS4}.$\quad\Box$

{\bf Теорема 1.2} {\sl  Пространства $l_D^{\,\theta}$ образуют шкалу компактно вложенных друг в друга пространств,
замкнутых относительно интерполяции (т.е. $[l_D^{\,0},\, l_D^{\,\theta}]_\tau =l_D^{\,\theta\tau}$ при всех
$\theta\ge0$, $\tau\in [0,1])$. При любом $\theta\geqslant 0$ оператор $T$ изоморфно отображает пространство
$W_2^\theta \ominus\{1\}$ на \ $l_D^{\,\theta}$.}

{\bf  Доказательство}
первого утверждения этой теоремы  полностью повторяет
доказательство Предложения 4 из работы
\cite{SS2}. Второе утверждение   доказано в Предложении 3 работы
\cite{SS4}.$\quad\Box$

Следующая теорема является ключевой для доказательства заключительной теоремы работы о равномерной устойчивости. В
частности, она говорит, что рассматриваемое отображение $F$ является слабо нелинейным, т.е. компактным возмущением
линейного отображения (при $\theta>0$).

 {\bf Теорема 1.3} {\sl При любом фиксированном $\theta \geqslant 0$
 оператор \ $F$\  отображает пространство $W^\theta_2$\
 в \  $l_D^{\,\theta}$\  и допускает представление вида
 $$
 F\,(\sigma) = T\ \sigma + \Phi(\sigma),
$$
где  \ $T$ \ --- линейный оператор, определенный в Теореме 1.1, а $\Phi$\ отображает пространство \ $W_2^\theta$
\ в $l^{\,\tau}_D$, \ где
$$
\tau=\begin{cases}2\theta,\quad& \text{если}\
0\leqslant\theta\leqslant1,\\ \theta+1,\quad& \text{если}\
1\leqslant\theta<\infty.\end{cases}
$$
Кроме того, отображение  \ $\Phi\, :\, W_2^\theta\to
l_D^{\,\tau}$\ является ограниченным в любом шаре,
т.e.
\begin{equation*}
\|\Phi(\sigma)\|_\tau\le C(R),\quad\text{если}\ \|\sigma\|_\theta\le R,
\end{equation*}
где постоянная \ $C$\ зависит только от радиуса шара
$R$.}

{\bf  Доказательство} этой теоремы проведено в работе
\cite{SS4}.$\quad\Box$

Теперь приведем  важные результаты об операторах
преобразования, которые хорошо известны в
классическом случае, а для сингулярных потенциалов
$q\in W^{-1}_2$  получены Гринивым и Микитюком
\cite{HM1} - \cite{HM2}.

{\bf Теорема 1.4.} (см. \cite{HM1}) {\sl Пусть функции $\sigma$ и $\tau$  вещественны и принадлежат пространству
$L_2[0,\pi]$. Пусть $L_D(\sigma)$  и $L_D(\tau)$  --- операторы Штурма--Лиувилля, порожденные дифференциальными
выражениями $-y''+\sigma'y$ и $-y''+\tau'y$ соответственно и краевыми условиями Дирихле $y(0)=y(\pi)=0$. Обозначим
через ${\rm s}_\sigma(x,\la)$ и ${\rm s}_\tau(x,\la)$ решения задач Коши $-y''+\sigma'y=\la y$ и $-y''+\tau'y=\la y$ с начальными условиями ${\rm s}(0,\la)=0$, ${\rm s}^{[1]}(0,\la)=\sqrt{\la}$. Тогда найдется функция $k_{\sigma,\tau}(x,t)$ с конечной нормой
$$
\|k\|:=\max\{\max\limits_{x\in[0,\pi]}\|k(x,\cdot)
\|_{L_2},\,\max\limits_{t\in[0,\pi]}\|k(\cdot,t)\|_{L_2}\},
$$
такая, что
$$
{\rm s}_\tau(x,\la)={\rm s}_\sigma(x,\la)+\int_0^xk_{\sigma,\tau}(x,t){\rm s}_\sigma(t,\la)dt.
$$}

{\bf Теорема 1.5} (см. \cite{HM2}) {\sl Пусть $\{\la_{k,\sigma}\}_1^\infty$ и $\{\la_{k,\tau}\}_1^\infty$ ---
собственные значения операторов $L_D(\sigma)$ и $L_D(\tau)$, а
$\{\varphi_{k,\sigma}={\rm s}_\sigma(x,\la_{k,\sigma})\}_1^\infty$ и $\{\varphi_{k,\tau}={\rm s}_\tau(x,\la_{k,\tau})\}_1^\infty$  --- соответствующие собственные функции. Положим $\al_{k,\sigma}=\|\varphi_{k,\sigma}\|^2$, $\al_{k,\tau}=\|\varphi_{k,\tau}\|^2$ и определим функцию
$$
F_{\sigma,\tau}(x,t)=
\suml_{k=1}^\infty\left(\frac1{\al_{k,\tau}}\varphi_{k,\tau}(x)\varphi_{k,\tau}(t)-
\frac1{\al_{k,\sigma}}\varphi_{k,\sigma}(x)\varphi_{k,\sigma}(t)\right).
$$
Тогда выполнено уравнение Гельфанда--Левитана--Марченко
\begin{equation}\label{GLM}
k_{\sigma,\tau}(x,t)+F_{\sigma,\tau}(x,t)+\intl_0^x
k_{\sigma,\tau}(x,s)F_{\sigma,\tau}(s,t)ds=0.
\end{equation}}

{\bf Теорема 1.6} (см. \cite{HM2}) {\sl Предположим, что функция $F_{\sigma,\tau}(x,t)$ непрерывна на квадрате, а функция $k_{\sigma,\tau}(x,t)$ непрерывна при $t<x$. Тогда
\begin{equation}\label{ktosigma}
\tau(x)-\sigma(x)=2k_{\sigma,\tau}(x,x)+C,
\end{equation}
где $C$ --- некоторая постоянная.}

Отметим, что в работах \cite{HM1} и \cite{HM2} последние три теоремы доказаны для случая $\sigma=0$. Переход к
произвольной функции $\sigma$ можно осуществить приемом из книги \cite[гл.1, \S1]{Le}.

\smallskip

\section{\bf Другие свойства отображения $F$.
Характеризация спектральных данных. Равномерная
устойчивость.
 }

Мы уже отмечали, что в рассматриваемой задаче функцию
$\sigma$ можно восстанавливать лишь с точностью до
константы. Поэтому далее удобнее работать не с
пространством $\sigma\in W^\theta_2\ominus\{1\}$, а с
фактор-пространством $\sigma\in W^\theta_2 / \{1\}$.
 Подразумеваем, что скалярное
произведение функций  $f,g\in  W^\theta_2 / \{1\}$
определено равенством $(f,g)_\theta =(f_0,
g_0)_\theta$,  где $f_0,g_0\in
W^\theta_2\ominus\{1\}$.

 Обозначим
через $\Gamma^{\,\theta}$ множество вещественных
функций
 $\sigma\in W^\theta_2 /\{1\}$, для
которых $\la_1(\sigma)\geqslant  1/2$, а через
$\mathcal B^{\,\theta}_\Gamma (R)$ --- пересечение
множества $\Gamma^{\,\theta}$ с замкнутым шаром
$\mathcal B^{\,\theta}_{\mathbb R}(R)$ радиуса $R$ в
пространстве $W^{\,\theta}_{2, \mathbb R}$. Здесь
число $1/2$ выбрано для определенности и простоты,
вместо $1/2$  может участвовать любое число $\eta>0$,
но тогда в \eqref{s1} и \eqref{sh} вместо
$s_2\geqslant 0$   нужно писать $s_2\geqslant
\sqrt\eta -1$.

 Если $\sigma\in\Gamma^{\,\theta}$, то собственные значения
оператора $L_D$ подчинены условиям $1/2
\leqslant \lambda_1 <\la_2<\dots$.
Для регуляризованных спектральных данных эти
неравенства эквивалентны следующим
\begin{equation}\label{s1}
s_2 \geqslant 0, \qquad s_{2k}-s_{2k+2}<1, \qquad
k=1,2,\dots .
\end{equation}
Условия неотрицательности всех нормировочных чисел
эквивалентны условиям
\begin{equation}\label{s2}
s_{2k-1}>-\pi/2,\qquad k=1,2,\dots.
\end{equation}
Последовательность $\{s_k\}_1^\infty \in l_2$,
поэтому  для любой вещественной функции $\sigma\in
\Gamma^{\,\theta}_D$ найдется число $h=h(\sigma)>0$,
такое, что
\begin{equation}\label{sh}
s_2 \geqslant 0, \qquad s_{2k}-s_{2k+2}\leqslant
1-h,\qquad s_{2k-1}\geqslant -\pi/2+h,
\qquad k=1,2,\dots .
\end{equation}

Фиксируем произвольные числа $r>0$ и $h \in (0,1)$.
Обозначим через $\Omega^{\,\theta}(r,h)$ совокупность
вещественных последовательностей $\{s_k\}_1^\infty$,
для которых выполнены неравенства
\eqref{sh} и которые лежат в замкнутом шаре радиуса
$r$ пространства $l_D^{\,\theta}$,  т.е.
$\|\{s_k\}\|_\theta
\leqslant r$  (здесь и далее подразумеваем, что $\|\cdot\|_\theta$
означает норму в пространстве  $l_D^{\,\theta}$).
Через $\Omega^{\,\theta}$ обозначим множество всех
вещественных последовательностей $\{s_k\}_1^\infty
\in l_D^{\,\theta}$, для которых справедливы неравенства
\eqref{s1} и \eqref{s2}.

Для доказательства  Теоремы 2.5  нам понадобится
следующий важный результат, который дает явное
описание  прообраза отображения $F$ при изменении
только одной из координат в пространстве
$l^{\,\theta}_D$. Похожие формулы для задачи
восстановления по одному спектру имеются в книге
\cite{PoTr}. Но доказательство нашего результата
проводится на другом пути.


{\bf Лемма 2.1} {\sl Пусть $\{\la_k\}$ и $\{\al_k\}$ --- собственные значения и нормировочные числа оператора $L_D$ с
функцией $\sigma\in W_2^\theta$, $\theta\ge0$. Тогда для любого фиксированного $n\ge1$ и для любого
$\xi\in(\la_{n-1}-\la_n,\la_{n+1}-\la_n)$ существует функция $\sigma(x,\xi)\in W_2^\theta$ такая, что оператор $L_D$,
построенный по этой функции имеет спектр $\{\la_k+\xi\delta_{kn}\}_1^\infty$ (здесь $\delta_{kn}$ --- символ Кронекера)
и нормировочные числа $\{\al_k\}_1^\infty$. Далее, для любого фиксированного $n\ge1$ и для любого
$\xi\in(-\al_n,+\infty)$ существует функция $\sigma(x,\xi)\in W_2^\theta$ такая, что оператор $L_D$, построенный по
этой функции имеет спектр $\{\la_k\}_1^\infty$ и нормировочные числа $\{\al_k+\xi\delta_{kn}\}_1^\infty$.}

{\bf Доказательство.} Для случая классических потенциалов эту лемму можно вывести из уравнения
Гельфанда--Левитана--Марченко с помощью известного в теории солитонов приема Йоста--Кона (см., например, \cite[гл.
2.5.7.]{Le}) В общем случае воспользуемся сформулированными в предыдущем параграфе теоремами 1.4 - 1.6.

В первом случае, когда меняется собственное значение
$\la_n\mapsto \la_n+\xi$ имеем
$$
F(x,t)=\frac1{\al_n}({\rm s}_\sigma(x,\la_n+\xi){\rm s}_\sigma(t,\la_n+\xi)-
{\rm s}_\sigma(x,\la_n){\rm s}_\sigma(t,\la_n).
$$
В этом случае уравнение Гельфанда--Левитана--Марченко
легко решается. Действительно, будем искать решение
--- функцию $k(x,t)$ в виде
$k(x,t)=c_1(x){\rm s}_\sigma(t,\la_n)+c_2(x){\rm s}_\sigma(t,\la_n+\xi)$, где $t<x$. Тогда уравнение \eqref{GLM} сводится к линейной системе
$$
\begin{cases}
c_1(x)-\frac1{\al_n}{\rm s}_\sigma(x,\la_n)-\frac1{\al_n}c_1(x)\intl_0^x{\rm s}^2_\sigma(t,\la_n)dt-
\frac1{\al_n}c_2(x)\intl_0^x{\rm s}_\sigma(t,\la_n){\rm s}_\sigma(t,\la_n+\xi)dt=0,\\
c_2(x)+\frac1{\al_n}{\rm s}_\sigma(x,\la_n+\xi)+\frac1{\al_n}c_1(x)\intl_0^x{\rm s}_\sigma(t,\la_n){\rm s}_\sigma(t,\la_n+\xi)dt+
\frac1{\al_n}c_2(x)\intl_0^x{\rm s}_\sigma^2(t,\la_n+\xi)dt=0.
\end{cases}
$$
Пользуясь формулой \eqref{ktosigma} потенциал $\sigma(x,\xi)$ можно написать в явном виде
\begin{equation}\label{varsigma_la}
\sigma(x,\xi)=\sigma(x)- 2\tfrac{d}{dx}\ln G(x,\xi),
\end{equation}
где
\begin{multline*}
G(x,\xi) =\left(1+\al_n^{-1}\int_0^x {\rm s}_\sigma^2(t,\la_n+\xi)dt\right)\left(1-\al_n^{-1}\int_0^x
{\rm s}_\sigma^2(t,\la_n)dt\right)\\
+\left(\al_n^{-1}\int_0^x{\rm s}_\sigma(t,\la_n+\xi){\rm s}_\sigma(t,\la_n)dt\right)^2.
\end{multline*}
Напомним, что здесь ${\rm s}_\sigma(x,\la)$  есть решение
уравнения $-y''+\sigma'y=\la y$ с начальными
условиями ${\rm s}_\sigma(0,\la)=0$, ${\rm s}_\sigma^{[1]}(0,\la)=\sqrt{\la}$.

При изменении одного нормировочного числа формула
имеет более простой вид:
\begin{equation}\label{varsigma_al}
\sigma(x,\xi)=\sigma(x)- 2\tfrac{d}{dx}\ln G(x,\xi),
\end{equation}
где
\begin{equation*}
G(x,\xi)=1+((\al_n+\xi)^{-1}-\al_n^{-1})\int_0^x{\rm s}_\sigma^2(t,\la_n)dt.
\end{equation*}
$\Box$

{\bf Лемма 2.2} {\sl При любом фиксированном $\theta \in[0,1/2)$\
отображение \ $F_D: \Gamma^\theta
\to \Omega^\theta$\ сюръективно.}

{\bf Доказательство.}
 Воспользуемся приемом
из  книги  Пошеля и Трубовица \cite{PoTr}. Согласно Теоремам 1.1 и 1.3 производная по Фреше отображения $F$ в точке
$\sigma =0$ совпадает с оператором $T$, который является изоморфизмом. В силу теоремы об обратном отображении для
 любого достаточно малого числа $\varepsilon >0$ найдется такое $\delta>0$, что образ шара $\|\sigma\|_\theta<\delta$ при отображении $F$
накрывает шар $\|s\|_\theta<\varepsilon$. При $\theta<1/2$ пространство $l_D^{\,\theta}$ совпадает с пространством
$l_2^{\,\theta}$.  Для данного $\bold s=\{s_k\}\in\Omega^\theta$ рассмотрим последовательность
$$
\bold s^n =\{0,0,\dots,0,s_n,s_{n+1},\dots\},
$$
выбрав число $n$ так, чтобы $\|\bold s^n
\|_\theta<\varepsilon$. Тогда найдется единственная
функция $\sigma_n\in W_2^\theta$, образ которой
$F(\sigma_n)$ совпадает с $\bold s^n$. Применив лемму
2.1 ($n-1$) раз, построим функцию
$\sigma\in\Gamma^{\,\theta}\subset W_{2,\mathbb
R}^\theta$, для которой $F \sigma =
\bold s$. Это и означает, что образ отображения $F$
содержит $\Omega^\theta$. Лемма доказана.$\quad\Box$

{\bf Лемма 2.3} {\sl Отображение \ $F_D: \Gamma^\theta \to \Omega^\theta$\
сюръективно при любом фиксированном $\theta \ge0$.}

{\bf Доказательство} Если при $\theta\in[0,1/2)$ сюръективность
уже доказана, то переход к произвольному $\theta\ge0$
осуществляется с помощью теоремы 1.3 дословным
повторением  доказательства теоремы 2.1 нашей работы
\cite{SS5}.$\quad\Box$

{\bf Лемма 2.4} {\sl При любых фиксированных
$\theta \ge0$\ отображение \ $F: \Gamma^\theta \to
\Omega^\theta$\ инъективно. }

{\bf Доказательство} Сошлемся на работу Гринива и Микитюка \cite{HM2}, где инъективность доказана для $\theta=0$, а
значит и для всех остальных $\theta\ge0$. Отметим также, что доказательство этой леммы можно провести методом работы
\cite{SS2} (см. Лемму 6) с учетом доказанной
ниже леммы 2.9.$\quad\Box$

{\bf Теорема 2.5} {\sl При любом фиксированном $\theta \geqslant 0$\ отображение \ $F: \Gamma^{\,\theta}
\to \Omega^{\,\theta}$\ есть биекция.  }

{\bf Доказательство.} Утверждение теоремы есть следствие лемм 2.3 и 2.4.$\quad\Box$

Обозначим через $\widehat \Omega^{\,\theta}$
множество последовательностей $\{s_k\}_1^\infty \in
l^{\,\theta}_D$, для которых числа $\lambda_k
=(s_{2k} +k)^2$ образуют строго возрастающую
последовательность, а все числа  $\alpha_k =
s_{2k-1}+\pi/2 $ положительны. Заметим, что если к
функции $\sigma$, которой определяется оператор $L_D$
, добавить функцию $cx$, то этот оператор перейдет в
$L_D +c$, т.е. его спектр сдвинется на $c$, а
нормировочные числа умножатся на коэффициенты
$\tfrac{\sqrt{\la_n+c}}{\sqrt{\la_n}}$. Положим
\begin{equation}\label{shift}
s_{2k}(c) = \sqrt{\lambda_k +c} -k,\qquad
s_{2k-1}(c)=\tfrac{\sqrt{\la_k+c}}{\sqrt{\la_k}}(s_{2k-1}+
\tfrac\pi2)-\tfrac\pi2.
\end{equation}
При этом если одно из $\la_n=0$, то, согласно
определению, нормировочное число
$\al_n(c)=\sqrt{c}\al_n$, а соответствующее число
$s_{2n-1}(c)=\sqrt{c}(s_{2n-1}+\pi/2)-\pi/2$.
Поскольку $cx
\in W^\theta_2$ при всех $\theta\geqslant 0$,  то $\{s_k(c)\}^\infty_1 \in l^{\,\theta}_D$ если и только если
$\{s_k(0)\}^\infty_1 \in l^{\,\theta}_D$. Следовательно, $\{s_k\}^\infty_1
\in \widehat\Omega^\theta$ если и только если найдется $c\geqslant0$,  такое, что  $\{s_k(c)\}^\infty_1 \in
\Omega^\theta$. Из сделанных замечаний следует

{\bf Теорема  2.6.} {\sl Отображение $F:
W^\theta_{2,\mathbb R}/\{1\} \to
\widehat\Omega^{\,\theta}$ есть биекция. Числа
$\{\lambda_k\}^\infty_1$ и $\{\alpha_k\}^\infty_1$
представляют спектр и нормировочные числа оператора
$L_D$  если и только если первая последовательность
является строго монотонной, вторая состоит из
положительных чисел
 и $\{s_k\}^\infty_1 \in
l_D^{\,\theta}$. }

Далее существенно  используются  аналитические свойства отображения $F$.  Мы предполагаем, что читатель знаком с
определением производных по Фреше и Гато для отображения $F: U\to H$, где $U$ --- открытое множество в $E$, а $E$  и
$H$  --- сепарабельные гильбертовы пространства. Для комплексных гильбертовых пространств производная по Фреше
естественно определяется в комплексном смысле. Отображение $F: U\to H$  называется аналитическим, если существует
комплексная Фреше производная в каждой точке $x\in U$. Фреше производную в точке $x$ далее обозначаем через $F'(x)$.
Естественным образом определяется понятие вещественного аналитического отображения, см., например,
\cite{PoTr}. Отображение $F: U\to H$ называется слабо
аналитическим, если в комплексном смысле
дифференцируемы по Гато координатные функции $(F(x),
e_k)$, где $\{e_k\}_1^\infty$  --- ортонормированный
базис пространства $H$. Известен результат
\cite{PoTr}, который значительно упрощает проверку
аналитичности отображения

{\bf Предложение 2.7}. {\sl Если $F: U\to H$ слабо
аналитическое отображение и локально  в каждой точке
$x\in U$  ограничено,  то $F$ ---  аналитическое
отображение.}

\vspace {0.2cm}

Далее мы говорим об отображениях замкнутых множеств (множеств вида $\Omega(r,h)$). Чтобы не делать дополнительных
объяснений, всюду считаем, что отображение $F:\, D\to H$ аналитично на $D$,  если найдется открытое множество $U$,
такое, что $U\supset D$ и $F:\, U\to H$ аналитично.

{\bf Теорема 2.8.} {\sl Пусть $\theta\geqslant 0$ и
$\sigma\in\Gamma^{\,\theta}$. Тогда найдется
комплексная окрестность $U\in W_2^\theta$ точки
$\sigma$, такая, что отображение $F: U
\to l_D^{\,\theta}$  является вещественно аналитическим.
  В этой окрестности отображение $\Phi =F-T: U \to l^{\,\tau}_D
$, где $\tau$ определено в Теореме 1.3, также
является вещественно аналитическим. Производная в
точке $\sigma\in U$ определяется равенством
\begin{equation} \label{F'D}
F'(\sigma) f  = \left\{ (\varphi_k(x),
\overline{f(x)})
\right\}_{k=1}^\infty,
\end{equation}
где
\begin{equation}\label{phi_k}
\varphi_{2k-1}(x) =
2\al_k\la_k\frac{d}{d\la}(z(x,\la)z'(x,\la))\vert_{\la
=\la_k},
\quad\varphi_{2k}(x)=-\frac{y'_k(x)y_k(x)}{\al_k\sqrt{\la_k}},
\qquad k=1,2, \dots .
\end{equation}
Здесь $f\in W^\theta_2$ --- функция, на которую действует оператор $ F'(\sigma): W_2^\theta\to l^{\,\theta}_D$,
$y_{n}=y(x,\la_n)$ --- собственные функции оператора $L_D$, нормированные условиями $y^{[1]}(0,\la_n)=\sqrt{\la_n}$, а
$z(x,\la)$ --- решение уравнения $-y''+\sigma'(x)y=\la y$ с начальным условием $z(\pi,\la)=0$, нормированное условием
$\int_0^\pi z^2(x,\la)dx=\frac1\la$. Утверждение об аналитичности (обычной) сохраняется, если условие
$\sigma\in\Gamma^{\,\theta}$ заменить условием $\sigma\in W^{\,\theta}_{2,\mathbb R}$ и потребовать, чтобы ноль не был
собственным значением оператора $L_D$.}

{\bf Доказательство.} В силу теоремы 1.3 и предложения 2.7
 достаточно доказать, что оператор, определенный в
\eqref{F'D} является производной отображения $F$ по Гато.
При этом половина формул уже доказана в нашей работе
\cite{SS2}. Действительно, если
$\sigma(t)=\sigma+tf$, где $t$ --- малый комплексный
параметр, то
$$
\dot{\la}_k(0)=-\frac{\left(y'_k(x)y_k(x),
\overline{f(x)}\right)}{(y^2_k,1)}.
$$
Таким образом, необходимо лишь найти значения
производных $\dot{\al}_n(0)$. Это требует
технической работы.

Обозначим через $\v(x,\la)$ решение уравнения
$-\v''+\sigma'\v=\la\v$ с граничными условиями
$\v(\pi,\la)=0$, $\v^{[1]}(0,\la)=1$. Такое решение
определено при всех комплексных $\la$ за исключением
собственных значений $\nu_n$ задачи Неймана--Дирихле.
Нам будет достаточно существования этого решения в
малых окрестностях $|\la-\la_n|<\eps_n$ точек
$\la_n$. Такие окрестности существуют, поскольку
числа $\la_n$ и $\nu_n$ различны. Через $\w(x,\la)$
мы обозначим другое решение уравнения
$-\w''+\sigma\w=\la\w$ с начальными условиями
$\w(0,\la)=1$, $\w^{[1]}(0,\la)=0$. Сразу же заметим,
что функции $\v$ и $\w$ линейно независимы, а их
Вронскиан равен $-1$. Положим
$\al(\la):=\la\int_0^\pi\v^2(x,\la)dx$ и заметим, что
в точках $\la_n$ значения этой функции совпадают с
нормировочными числами $\al_n$. Далее нам потребуется
выражение для решения неоднородного уравнения
$-y''+\sigma'y=\la y+f$ с правой частью $f\in
L_1[0,\pi]$ и краевыми условиями
$y(\pi,\la)=y^{[1]}(0,\la)=0$. Легко проверить, что
\begin{equation}\label{neodnorod}
y(x,\la)=-\v(x,\la)\intl_0^x\w(\xi,\la)f(\xi)d\xi-\w(x,\la)\intl_x^\pi\v(\xi,\la)f(\xi)d\xi.
\end{equation}
В частности мы можем найти выражение для функции $\v_\la(x,\la)$ --- производной функции $\v(x,\la)$ по $\la$ (во
избежание недоразумений далее штрихом обозначаем производную по $x$, точкой --- производную по $t$ и индексом $\la$
--- производную по $\la$). Дифференцируя уравнения для функции $\v$ по переменной $\la$ придем к равенству
$$
-\v_\la''+\sigma'\v_\la=\la\v_\la+\v,\qquad \v_\la(\pi,\la)=\v_\la^{[1]}(0,\la)=0.
$$
Тогда имеем
$$
\v_\la(x,\la)=-\v(x,\la)\intl_0^x\w(\xi,\la)\v(\xi,\la)d\xi-\w(x,\la)\intl_x^\pi\v^2(\xi,\la)d\xi.
$$
Равенство \eqref{neodnorod} позволяет также получить выражение для функции $\dot{\v}(x,\la)$. Действительно,
дифференцируя уравнение $-\v''+\sigma'\v=\la_n\v$ по переменной $t$ в точке $t=0$ и учитывая, что $\sigma=\sigma+tf$, а
$\la_n$ также есть функция переменной $t$, получим
$$
-\dot\v''+\sigma'\dot\v=\la_n\dot\v+(\dot\la_n-f')\v,\qquad \dot\v(\pi,\la_n)=\dot\v^{[1]}(0,\la_n)=0.
$$
Тогда
$$
\dot\v(x,\la_n)=-\v(x,\la_n)\intl_0^x\w(\xi,\la_n)(\dot\la_n-f'(\xi))\v(\xi,\la_n)d\xi-\w(x,\la_n)\intl_x^\pi
\v^2(\xi,\la_n)(\dot\la_n-f'(\xi))d\xi.
$$
Дифференцируя по $t$ выражение $\al_n=\al(\la_n)$, получим
\begin{equation}\label{aldot1}
\dot\al_n=\dot\la_n\intl_0^\pi\v^2(x,\la_n)dx+2\la_n\intl_0^\pi\v(x,\la_n)\dot\v(x,\la_n)dx.
\end{equation}
Преобразуем последний интеграл
\begin{gather*}
\intl_0^\pi\v(x,\la_n)\dot\v(x,\la_n)dx=-\intl_0^\pi
\v^2(x,\la_n)\intl_0^x\w(\xi,\la_n)\v(\xi,\la_n)(\dot\la_n-f'(\xi))d\xi-\\-\intl_0^\pi\v(x,\la_n)\w(x,\la_n)\intl_x^\pi
\v^2(\xi,\la_n)(\dot\la_n-f'(\xi))d\xi=
\end{gather*}
\begin{gather*}=-\intl_0^\pi(\dot\la_n-f'(\xi))\v(\xi,\la_n)\left[\w(\xi,\la_n)\intl_\xi^\pi
\v^2(x,\la_n)dx+\v(\xi,\la_n)\intl_0^\xi
\v(x,\la_n)\w(x,\la_n)dx\right]d\xi=\\=\intl_0^\pi(\dot\la_n-f'(\xi))\v(\xi,\la_n)\v_\la(\xi,\la_n)d\xi.
\end{gather*}
Подставляя полученное выражение в формулу \eqref{aldot1} и учитывая, что
$$
\dot\la_n=\frac{\intl_0^\pi\v^2(\xi,\la_n)f'(\xi)d\xi}{\intl_0^\pi\v^2(x,\la_n)dx},
$$
получим
\begin{equation}\label{aldot2}
\dot\al_n=\left[1+\frac{2\la_n\intl_0^\pi\v(x,\la_n)\v_\la(x,\la_n)dx}{\intl_0^\pi\v^2(x,\la_n)dx}\right]\intl_0^\pi\v^2(\xi)f'(\xi)-
2\la_n\intl_0^\pi\v(\xi,\la_n)\v_\la(\xi,\la_n)f'(\xi)d\xi.
\end{equation}
Заметим теперь, что
$$
\al_\la(\la_n)=\intl_0^\pi\v^2(x,\la_n)dx+2\la_n\intl_0^\pi\v(x,\la_n)\v_\la(x,\la_n)dx,
$$
а тогда \eqref{aldot2} можно преобразовать к виду
\begin{equation}\label{aldot3}
\dot\al_n=-\la_n\al_n\frac{d}{d\la}\left(\frac{\intl_0^\pi\v^2(\xi,\la_n)f'(\xi)d\xi}{\al(\la)}\right).
\end{equation}
Поскольку $\int_0^\pi\v^2(x,\la)dx=\frac{\al}{\la}$, то $\v(x,\la)=z(x,\la)\sqrt{\al(\la)}$, а значит окончательно
имеем
$$
\dot\al_n=-\la_n\al_n\frac{d}{d\la}\left(\intl_0^\pi z^2(\xi,\la_n)f'(\xi)d\xi\right)=
\la_n\al_n\frac{d}{d\la}\left(\intl_0^\pi2z(\xi,\la_n)z'(\xi,\la_n)f(\xi)d\xi\right).
$$
Теорема доказана.$\quad\Box$

{\bf Лемма 2.9}. {\sl В обозначениях леммы 2.8
 система функций $\{\varphi_k\}_1^\infty$, где
\begin{equation}\label{phiD}
\varphi_{2k-1}(x) =
2\al_k\la_k\frac{d}{d\la}(z(x,\la_k)z'(x,\la_k)),\qquad\varphi_{2k}(x)=-\frac{y'_k(x)y_k(x)}{\al_k\sqrt{\la_k}},
\qquad k=1,2, \dots,
\end{equation}
является базисом Рисса в пространстве $L_2(0, \pi)/\{1\}$, квадратично близким к базису $\left\{(x-\pi)\cos(2kx),\,-\frac
1\pi \sin (2kx)\right\}_1^\infty$. Биортогональная система к $\{\varphi_k(x)\}_1^\infty$ имеет вид
$\{\psi_k(x)\}_1^\infty$, где
\begin{equation}\label{psiD}
\psi_{2k-1}(x)=\frac2{\al_k^2}y_k^2(x),\qquad\psi_{2k}(x)=-\frac{2\sqrt{\la_k}}{\al_k}\frac{d}{d\la}\left(y^2(x,\la_k)\right).
\end{equation}}

{\bf Доказательство.} Вначале мы докажем соотношения биортогональности. Рассмотрим функцию
$$
\delta_m(\la):=\intl_0^\pi y^2(x,\la)(y^2_m(x))'dx.
$$
Проведя интегрирование по частям, получим $\delta_m(\la)=-\int_0^\pi(y^2(x,\la))'y^2_m(x)dx$. Тогда
\begin{multline*}
2\delta_m(\la)=\intl_0^\pi\left(y^2(x,\la)(y^2_m(x))'-(y^2(x,\la))'y_m(x)\right)dx=\\
=2\intl_0^\pi y(x,\la)y_m(x)\left(y(x,\la)y'_m(x)-y'(x,\la)y_m(x)\right)dx.
\end{multline*}
При $\la=\la_m$ имеем $\delta_m(\la)=0$ и далее считаем $\la\ne\la_m$. Поскольку
$\left(y(x,\la)y'_m(x)-y'(x,\la)y_m(x)\right)'=(\la-\la_m)y(x,\la)y_m(x)dx$, то
\begin{multline*}
2\delta_m(\la)=\frac1{\la-\la_m}\intl_0^\pi\left(y(x,\la)y'_m(x)-y'(x,\la)y_m(x)\right)'\left(y(x,\la)y'_m(x)-y'(x,\la)y_m(x)\right)dx=\\
=\frac1{\la-\la_m}\left[y(x,\la)y'_m(x)-y'(x,\la)y_m(x)\right]^2\vert_0^\pi=\frac{y'_m(\pi)}{\la-\la_m}y^2(\pi,\la).
\end{multline*}
Последнее выражение обращается в ноль вместе со своей первой производной в точке $\la=\la_n\ne\la_m$. Таким образом мы
доказали, что $(\varphi_{2n},\psi_{2m-1})=0$ при любых $n$ и $m$, а $(\varphi_{2n},\psi_{2m})=0$ при любых $n\ne m$.

Теперь рассмотрим функцию
$$
\delta(\la,\mu):=\intl_0^\pi y^2(x,\la)(z^2(x,\mu))'dx.
$$
Проводя аналогичные рассуждения при $\mu\ne\la$ придем к равенству
$$
2\delta(\la,\mu)=\frac1{\la-\mu}\left[(y(\pi,\la)z'(\pi,\mu))^2-(y'(0,\la)z(0,\mu))^2\right].
$$
Таким образом, $\delta(\la,\mu)=\varkappa_1(\la,\mu)(\la-\la_n)^2+\varkappa_2(\la,\mu)(\mu-\la_m)^2$, где функции
$\varkappa_j(\la,\mu)$ аналитичны в окрестности точки $\la=\la_n$, $\mu=\la_m$. из этого представления следует, что
частные производные $\tfrac{d}{d\la}$ и $\tfrac{d^2}{d\la d\mu}$ функции $\delta(\la,\mu)$ обращаются в ноль в этой
точке . Это доказывает соотношения $(\varphi_{2n-1},\psi_{2m})=(\varphi_{2n-1},\psi_{2m-1})=0$ при всех $n\ne m$.

Доказательство равенства $(\varphi_{2n-1},\psi_{2n})=0$ требует других рассуждений, поскольку функция
$\delta(\la,\mu)$, используемая выше не определена при $\la=\mu$. Запишем
$$
(\varphi_{2n-1},\psi_{2n})=-8\la_n\intl_0^\pi y(x,\la_n)y_\la(x,\la_n)(z(x,\la_n)z_\la(x,\la_n))'dx.
$$
Вновь интегрируя по частям, получим
\begin{multline*}
(\varphi_{2n-1},\psi_{2n})=4\la_n\intl_0^\pi\left[(y(x,\la_n)y_\la(x,\la_n))'z(x,\la_n)z_\la(x,\la_n)-\right.\\
\left.-y(x,\la_n)y_\la(x,\la_n)(z(x,\la_n)z_\la(x,\la_n))'\right]dx.
\end{multline*}
При $\la=\la_n$ функции $y(x,\la)$ и $z(x,\la)$ линейно зависимы, а значит
$$
(\varphi_{2n-1},\psi_{2n})=4\la_n\intl_0^\pi
y(x,\la_n)z(x,\la_n)(y_\la(x,\la_n)z'_\la(x,\la_n)-y'_\la(x,\la_n)z_\la(x,\la_n))dx,
$$
где индекс $\la$ означает дифференцирование по $\la$, а штрих --- дифференцирование по $x$. Заметим теперь, что
$$
(y(x,\la)z'(x,\mu)-y'(x,\la)z(x,\mu))'=(\la-\mu)y(x,\la)z(x,\mu)
$$
Дифференцируя по $\la$ и $\mu$ последнее равенство и устремляя $\la\to\la_n$, $\mu\to\la_n$ получим
$(\varphi_{2n-1},\psi_{2n})=0$.

Докажем наконец, что $(\varphi_n,\psi_n)=1$. Это равенство при четных $n$ принимает вид
$$
\frac2{\al_n^2}\intl_0^\pi(y^2(x,\la_n))_\la(y^2(x,\la_n))'dx=1.
$$
Вновь введем функцию
$$
\delta(\la)=\frac2{\al_n^2}\intl_0^\pi y^2(x,\la)(y^2(x,\la_n))'dx.
$$
Повторяя рассуждения, приведенные выше, получим
$$
\delta(\la)=\frac1{\al_n^2(\la-\la_n)}y^2(\pi,\la)(y_n'(\pi))^2.
$$
Тогда
$$
\delta'(\la_n)=\frac1{\al_n^2}y^2_\la(\pi,\la_n)(y_n'(\pi))^2=1.
$$
Здесь мы воспользовались равенством $\al_n=y_\la(\pi,\la_n)y'(\pi,\la_n)$, доказанным в работе \cite{SS4}.

Равенство $(\varphi_n,\psi_n)=1$ при нечетных $n$ доказывается аналогично после применения формулы
$\la_nz_\la(0,\la_n)=-\sqrt{\al_n}$.

Итак, мы доказали биортогональность систем $\{\varphi_k\}_1^\infty$ и $\{\psi_k\}_1^\infty$. Теперь достаточно доказать
квадратичную близость одной из этих систем к базису Рисса --- квадратичная близость другой системы (к биортогональному
базису) будет следовать тогда из общих теорем. Нам будет удобно провести доказательство для системы
$\{\psi_k\}_1^\infty$.

Асимптотические формулы $y_k(x)=\sin(kx)+\chi_k(x)$, где $\sum_{k=1}^\infty\|\chi_k\|_C^2<\infty$ доказаны в работе
\cite{SS2003}. В работе \cite{SS3} доказаны асимптотические соотношения для нормировочных чисел
$\al_k=\pi/2+\rho_k$, где $\sum_{k=1}^\infty|\rho_k|^2<\infty$. Отсюда следует, что
$$
\al_k^2y^2_k=\frac{\pi^2}4\sin^2(kx)+\chi_k,\qquad\text{где} \suml_{k=1}^\infty\|\chi_k\|_C^2<\infty.
$$

Нам, однако, нужны еще асимптотические выражения для функций $y_\la(x,\la_k)$. Обозначим через $s(x,\la)$ и $c(x,\la)$
решения уравнения $-y''+\sigma'(x)y=\la y$  с начальными условиями $s(0,\la)=c^{[1]}(0,\la)=0$, $c(0,\la)=1$,
$s^{[1]}(0,\la)=\sqrt{\la}$. Из результатов работы \cite{SS2003} следует, что
$$
s(x,\la)=\sin(\sqrt{\la}x)+\rho_1(x,\la),\qquad c(x,\la)=\cos\sqrt{\la}x+\rho_2(x,\la),
$$
где $\{\|\rho_j(x,\la_k)\|_C\}_{k=1}^\infty\in l_2$. Воспользуемся методом вариации постоянных и представим функцию
$y_\la(x,\la)$ --- решение уравнения $-y_\la''+\sigma'y_\la=\la y_\la+y$ с начальными условиями $y_\la(0,\la)=0$,
$y^{[1]}_\la(0,\la)=1/(2\sqrt{\la})$ в виде
$$
y_\la(x,\la)=\frac1{2\la}s(x,\la)+c(x,\la)\intl_0^x s^2(\xi,\la)d\xi-s(x,\la)\intl_0^xs(\xi,\la)c(\xi,\la)d\xi.
$$
Подставляя сюда асимптотические соотношения для функций $s(x,\la_n)$ и $c(x,\la_n)$ и учитывая, что
$\{\sqrt{\la_k}-k\}\in l_2$ придем к равенству
$$
y_\la(x,\la_k)=\frac{x\cos(kx)}{4k}+\frac1k\chi_k(x),\qquad\text{где}\ \suml_{k=1}^\infty\|\chi_k\|^2_C<\infty.
$$
Отсюда следуют асимптотические формулы
$$
-\frac{4\sqrt{\la_k
}}{\al_k}y(x,\la_k)y_\la(x,\la_k)=\frac2\pi x\sin(2kx)+\chi_k(x),\qquad\text{где}\
\suml_{k=1}^\infty\|\chi_k\|^2_C<\infty.
$$

Итак, мы доказали квадратичную близость системы $\{\psi_k\}_1^\infty$  к системе
$$
\left\{\frac{\pi^2}4\sin^2(kx),\,\frac2\pi x\sin(2kx)\right\}_1^\infty.
$$
В работе \cite{SS4} доказана базисность Рисса
последней системы (точнее в работе доказана
базисность Рисса биортогональной системы, что
равносильно). Лемма доказана.$\quad\Box$

{\bf Теорема 2.10.} {\sl Пусть $\theta\geqslant 0$.
Для каждой точки $\bold y_0
\in \Omega^\theta =F(\Gamma^\theta)$
 существует ее комплексная окрестность $U(\bold y_0)$,
в которой определено обратное отображение
$F^{-1}(\bold y)$ и в которой это отображение имеет
комплексную Фреше производную. Эта производная имеет
вид
\begin{equation}\label{F-1}
\left( F^{-1}\right)'(\bold y)= (F')^{-1}(\bold y) =
\sum^\infty_{k=1} s_k\psi_k(x), \qquad \bold y= (s_1, s_2, \dots ),
\end{equation}
где $\{\psi_k(x)\}_1^\infty$ ---  биортогональная
система из леммы 2.9. }

{\bf Доказательство} этой теоремы получается
дословным повторением теоремы 2.6 из нашей работы
\cite{SS5}.$\quad\Box$

Отметим, что из теорем 2.8 и 2.10  сразу получаются
локальные оценки разности потенциалов через разность
спектральных данных и наоборот. Во введении работы
\cite{SS5} отмечено,  что для классического случая
$\theta =1\ (q\in L_2)$
выполнено много работ на эту тему различными методами
и в разной форме. Однако, изучались отображения $q
\to \{\text{спектральные данные}\}$, мы же изучаем
отображение $\int q(t)\, dt = \sigma
\to \{\text{спектральные данные}\}$,
поэтому возникающие у нас системы и формулы имеют
другой вид.

Далее покажем, что при $\theta >0$  с помощью Теоремы
1.3 можно получить существенно более сильный
результат, избегая технической работы с системами
функций.  Доказательства носят общий характер и
получаются точно также, как в работе \cite{SS5}.

\vspace{0.2cm}

{\bf Лемма 2.11.} {\sl Фиксируем $\theta >0$.
Пусть $R$ произвольное положительное число и
$\mathcal B^{\,\theta}_\Gamma (R) = \Gamma
\cap\mathcal B^{\,\theta}_\mathbb R (R)$.
 Тогда найдутся положительные числа $r=r(R), h=h(R)$,
такие, что
$$
F (\mathcal B^{\,\theta}_\Gamma (R) ) \subset
\Omega^\theta (r,h).
$$
}

{\bf Доказательство} проводится также, как в лемме
2.7 работы \cite{SS5}.$\quad\Box$

{\bf Лемма 2.12.} {\sl Пусть $\theta >0$. Справедливо обратное утверждение к Лемме~2.11: для любых чисел $r$ и $h$
найдется число $R>0$,  такое, что
$$
F^{-1}(\,\Omega^{\,\theta}(r,h)) \subset \mathcal
B^{\,\theta}_\Gamma (R).
$$
Справедливо представление
$$
F^{-1} = T^{-1} +\Psi, \qquad \Psi: \Omega^{\,\theta}
\to W^\tau_2,
$$
где число $\tau$  определено в Теореме~1.3.
Отображение $\Psi: \Omega^\theta \to W^\tau_2,$
аналитично, причем
\begin{equation}\label{Psi}
\|\Psi\bold y\|_\tau \leqslant C\|\bold y\|_\theta
\qquad\text{для всех} \ \, \bold y \in\Omega^{\,\theta}(r,h),
\end{equation}
где постоянная $C$  зависит только от $r$  и $h$. }

{\bf  Доказательство}
проводится также, как в лемме 2.8 работы \cite{SS5}.$\quad\Box$

{\bf Лемма 2.13.} {\sl Пусть $\theta>0$. При любом $R>0$  справедлива
оценка
\begin{equation}\label{direct}
\|F'(\sigma)\|_\theta \leqslant C, \qquad \text{для всех}\ \ \sigma \in
\mathcal B_\Gamma^{\,\theta} (R),
\end{equation}
где постоянная $C$ зависит от $R$,  но не зависит от
$\sigma$. }

{\bf Доказательство}
проводится также, как в лемме 2.10 работы
\cite{SS5}.$\quad\Box$

{\bf Лемма 2.14.} {\sl Пусть $\theta >0$. При любых
$r>0,\ h\in (0,\, 1/2)$  для обратного отображения
справедлива оценка

\begin{equation}\label{inverse}
\|(F^{-1})'(\bold y )\| \leqslant C, \qquad \text{для всех}
\ \, \bold y \in \Omega^\theta(r,h),
\end{equation}
где постоянная $C$ зависит от $r$ и $h$,  но не
зависит от $\bold y$. }

{\bf Доказательство}
проводится также, как в лемме 2.11 работы
\cite{SS5}.$\quad\Box$

{\bf Теорема~2.15.} {\sl  Фиксируем $\theta >0$.
 Пусть последовательности
$\bold y, \bold y_1$  регуляризованных спектральных данных лежат в множестве $\Omega^{\,\theta}(r,h)$. Тогда  прообразы
$\sigma = F^{-1}\bold y,\ \,\sigma_1 = F^{-1}\bold y_1$  лежат в множестве $\mathcal B_\Gamma^{\,\theta} (R)$ и
справедливы оценки
\begin{equation}\label{a}
C_1\|\bold y -\bold y_1\|_\theta \leqslant
\|\sigma -\sigma_1\|_\theta
\leqslant C_2\|\bold y -\bold y_1\|_\theta,
\end{equation}
где число $R$ и постоянные $C_1, C_2$ зависят только
$r$ и $h$. Число $R$   и постоянные $C_2, C^{-1}_1$
увеличиваются  при $r\to\infty$  или $h\to 0$.
Обратно, если $\sigma,
\sigma_1$ лежат в шаре $\mathcal B_\mathbb
R^{\,\theta}(R)$, то последовательности  $\bold y,
\bold y_1$ регуляризованных спектральных данных этих функций лежат в множестве $\Omega^\theta(r,h)$  и справедливы оценки

\begin{equation}\label{b}
C_1 \|\sigma -\sigma_1\|_\theta \leqslant\|\bold y
-\bold y_1\|_\theta \leqslant C_2 \|\sigma
-\sigma_1\|_\theta .
\end{equation}
Здесь числа $r>0, \ \, h\in (0, 1/2)$  и постоянные $C_1$ и $C_2$ зависят только от $R$. Числа $ r,\, h^{-1}, C_2 $ и
$C^{-1}_1$  увеличиваются при $R\to\infty$. }

{\it Доказательство}
проводится также, как в теореме 2.12 работы
\cite{SS5}.$\quad\Box$

Множества  $\mathcal B_\Gamma^\theta(R)$ в Теореме
2.15 можно заменить обычными шарами $\mathcal
B^\theta_{\mathbb R}(R)$,  но тогда регуляризованные
спектральные данные нужно определить формулой
\eqref{shift},  где постоянная $c$  такова, что для
всех $\sigma \in \mathcal B^{\,\theta}_{\mathbb R}(R)
$ выполнена оценка $c\geqslant -
\lambda_1(\sigma)-1$. В силу теоремы 1.3 такая
постоянная, зависящая только от $R$ существует. Это замечание  вытекает из того, при добавлении к $\sigma$    функции
$cx$ спектр оператора $L_D$ сдвигается на $c$, а разность функций $\sigma ,\,
\sigma_1 \in
\mathcal B^{\,\theta}_{\mathbb R}(R)$ совпадает с разностью
функций $\sigma +cx, \ \sigma_1 +cx\in
\mathcal B_\Gamma^{\,\theta} (R)$.

\bigskip
\medskip
 \address{A.М.Савчук,
 МГУ имени М.В.Ломоносова,
механико-математический ф-т, Ленинские Горы, Москва,
119992.

\email{artem\_savchuk@mail.ru}}

\medskip

\address{A.А.Шкаликов, МГУ имени М.В.Ломоносова,
механико-математический ф-т, Ленинские Горы, Москва,
119992.

\email{ashkalikov@yahoo.com}
}

\end{document}